\documentclass[12pt]{article}
\usepackage{geometry,graphicx,amssymb,amsfonts,amsmath,bm,url}
\def\no{\noindent}

\newtheorem{theo}{Theorem}

\def\pmatrix{\left(\begin{array}}
\def\endpmatrix{\end{array}\right)}
\def\RR{{\mathbb{R}}}
\def\CC{{\mathbb{C}}}
\def\DD{{\mathbb{D}}}
\def\bfb{{\bm{b}}}
\def\bfc{{\bm{c}}}
\def\bfe{{\bm{e}}}

\def\bfv{{\bm{v}}}
\def\bfw{{\bm{w}}}
\def\bfy{{\bm{y}}}

\def\PP{{\cal P}}

\def\bb{{\beta}}

\def\dd{{\delta}}
\def\hDelta{\hat{\Delta}}
\def\hc{\hat{c}}
\def\hbfc{\hat{\bfc}}
\def\hPP{\hat{\PP}}
\def\hy{\hat{y}}
\def\hbfy{\hat{\bfy}}

\def\hG{\hat{G}}
\def\hL{\hat{L}}
\def\hU{\hat{U}}

\title{Efficient implementation of Radau collocation methods}

\author{Luigi Brugnano\footnote{Dipartimento di Matematica ``U.\,Dini'', Universit\`a di Firenze, Italy
({\tt luigi.brugnano@unifi.it})}
\and Felice Iavernaro\footnote{Dipartimento di Matematica, Universit\`a di Bari, Italy
({\tt felix@dm.uniba.it})}
\and Cecilia Magherini\footnote{Dipartimento di Matematica, Universit\`a di Pisa, Italy
({\tt magherini@dm.unipi.it})}
}

\date{December 29, 2012}
\begin{document}
\maketitle

\begin{abstract}
In this paper we define an efficient implementation of Runge-Kutta methods of Radau IIA type, which are commonly used
when solving stiff ODE-IVPs problems. The proposed implementation relies on an alternative low-rank formulation of the
 methods, for which a splitting procedure is easily defined. The linear convergence analysis of this splitting procedure exhibits excellent properties, which are confirmed by its performance on a few numerical tests.

\bigskip
\no{\bf Keywords:} Radau IIA collocation methods; W-transform; Implicit Runge-Kutta methods;
Singly Implicit Runge-Kutta methods; Splitting; Hamiltonian BVMs.

\bigskip
\no{\bf MSC (2010):} 65L04, 65L05, 65L06, 65L99.
\end{abstract}

\section{Introduction}
The efficient numerical solution of implicit Runge-Kutta methods has been the subject of many investigations in the last
 decades, starting from the seminal paper of Butcher \cite{Bu76,Bu79} (see also \cite{Bi77}). An $s$-stage R-K method
  applied to the initial value problem
\begin{equation}\label{ivp}
y' = f(y), \qquad y(t_0)=y_0\in\RR^m,
\end{equation}
yields a nonlinear system of dimension $sm$ which takes the form
\begin{equation}\label{R-K}
G(\bfy)\equiv \bfy - \bfe\otimes y_0 - h A \otimes I\, f(\bfy) = 0,
\end{equation}
where
\begin{equation}
\label{vectors}
\bfe=\pmatrix{c} 1 \\ \vdots \\1\endpmatrix\in\RR^s, \qquad \bfy=\pmatrix{c} y_1\\
\vdots\\ y_s\endpmatrix, \qquad f(\bfy)=\pmatrix{c} f(y_1)\\
\vdots\\ f(y_s)\endpmatrix,
\end{equation}
$y_1,\dots,y_s$ being the internal stages. It is common to solve \eqref{R-K} by a simplified Newton iteration, namely, for
 $k=0,1,\dots$,
\begin{equation}
\label{simnewt}
\begin{array}{l}
\left(I -hA\otimes J\right)\Delta^{(k)} = -G(\bfy^{(k)}), \\
\bfy^{(k+1)}=\bfy^{(k)}+\Delta^{(k)},
\end{array}
\end{equation}
where $J$ is the Jacobian of $f$ evaluated at some intermediate point and $\bfy^{(0)}$ an initial approximation of the stage
 vector, for instance $J=\frac{\partial f}{\partial y}(y_{0})$ and $\bfy^{(0)}=\bfe\otimes y_{0}$.
To reduce the computational efforts associated with the solution of \eqref{simnewt}, a suitable linear change of variables on
the $s$ stages of the method is often introduced with the goal of simplifying the structure of the system itself. This is
tantamount to performing a  similarity transformation, commonly referred to as \textit{Butcher transformation}, that puts the
coefficient matrix $A$ of the R-K method in a simpler form, i.e. a diagonal or triangular matrix. Let $B=TAT^{-1}$
such a transformation.
System \eqref{simnewt} becomes
\begin{equation}
\label{newt1}
\left( I - h(B\otimes J)\right)(T\otimes I) \Delta^{(k)} = -(T\otimes I)G(\bfy^{(k)}),
\end{equation}
with the obvious advantage that the costs associated with the $LU$ factorizations decrease from $O(s^{3}m^{3})$ to
$O(sm^{3})$ {\em flops}.\footnote{One {\em flop} is an elementary {\em fl}oating-point {\em op}eration.}
In particular, if $A$ has a one-point spectrum one only needs a single $LU$ decomposition and the cost further reduces to
 $O(m^3)$ {\em flops} \cite{Bu78}. However, for many fully implicit methods of interest, the matrix $A$ possesses complex
conjugate pairs  of eigenvalues which will appear as diagonal entries in the matrix $B$. In such a case, it is
computationally more advantageous to allow $B$ to be block-diagonal, with each $2\times 2$ diagonal block corresponding
 to a complex conjugate pair of eigenvalues of $A$. Each subsystem of dimension $2m$ is then turned into an \
 $m$-dimensional complex system.  This is the standard procedure used in the codes RADAU5 \cite{HaWa91,testset} and
  RADAU \cite{HaWa99,testset}, the former a variable-step fixed-order code, and the latter a variable-order variant,  both
   based upon Radau-IIA formulae (of orders $5$, $9$, and $13$).

Subsequent attempts to derive implicit high-order methods, for which the discrete problem to be solved can be cast in a
simplified form, have been made, e.g., in \cite{Al77,Ca83}. This line of investigation has been further refined in later papers
 (see, e.g., \cite{CoBu83,Co91}).
Sometimes, the formulation of the discrete problem has been suitably modified, in order to induce a corresponding
``natural splitting'' procedure, as is done, e.g., in \cite{Br00,BrMa02,BrMa04} (see also \cite{BrMa07,BrMaMu06}).

A different approach to the problem is that of considering suitable splitting procedures for solving the generated discrete
problems \cite{AmBr97,Seve06,HoSo91,HoSw97,HoSw97a,IaMa98,Sw97}. A particularly interesting splitting scheme, first
introduced in \cite{HoSw97}, is that induced by the Crout factorization of the coefficient matrix $A$, namely $A=LU$, with
$L$ lower triangular and $U$ upper triangular with unit diagonal entries. After observing that, for many remarkable R-K
methods, the lower triangular part of $A$ is \textit{dominant}, in \cite{HoSw97} the authors suggest to replace the matrix
$A$ in \eqref{simnewt} with the matrix, $L$ thus obtaining the scheme
\begin{equation}
\label{simnewt1}
\begin{array}{l}
\left(I -hL\otimes J\right)\Delta^{(k)} = -G(\bfy^{(k)}), \\
\bfy^{(k+1)}=\bfy^{(k)}+\Delta^{(k)}.
\end{array}
\end{equation}
Compared to \eqref{simnewt} this scheme only requires the sequential solution of $s$ subsystems of dimension $m$ and
 therefore a global cost of $O(sm^{3})$ elementary operations. Moreover, the $s$ $LU$ factorizations of the matrices
 $I-l_{ii}J$~ ($l_{ii}$ being the $i$th diagonal entry of $L$), and the evaluations of the components of $G(\bfy^{(k)})$ may be
 done in parallel. This is why  the corresponding methods have been named \textit{parallel triangularly implicit iterated
 R-K methods} (PTIRK).

On the other hand, if the original modified Newton process \eqref{simnewt} converges in one iterate on linear problems, the
 same no longer holds true for \eqref{simnewt1}, due to the approximation $A\simeq L$. Applying the method to the linear
  test equation $y'=\lambda y$ yields the following estimation for the error $e^{(k)} = y^{(k)}-y$:
\begin{equation}
\label{linerr}
e^{(k+1)}=M(q)e^{(k)}, \qquad M(q)=q(I-qL)^{-1}(A-L),
\end{equation}
with $q=h\lambda$. Matrix $M(q)$ is referred to as the \textit{amplification matrix} associated with the method and its
 properties influence the rate of convergence of the scheme \eqref{simnewt1} according to a first order convergence
  analysis (see Section \ref{four}).

In this paper we wish to combine both the approaches described above and epitomized at formulae \eqref{newt1} and
 \eqref{simnewt1}, to derive an efficient implementation of Radau IIA methods on sequential computers. In fact, the above
 discussion begs the following question: is it possible to perform a change of variables of the stage vector such that,
 for the new system \eqref{newt1}, the matrix $B$ admits a $LU$ factorization with constant diagonal entries?
 In the affirmative, a single $LU$ factorization would be needed to solve \eqref{simnewt1},
 with a cost of only $O(m^{3})$ {\em flops}. A first positive answer in this direction has been given in
 \cite{AmBr97} for general R-K methods. Later on, in \cite{Sw97}, an optimal splitting of the form \eqref{simnewt1}
 has been devised for the Radau IIA method of order three (two stages), with $l_{11}=l_{22}$.

In this paper, we follow a different route, which relies on a {\em low-rank} formulation of Radau IIA collocation methods.
Low-rank R-K methods have been recently introduced in a series of papers in the context of numerical geometric
integration  \cite{BrIaTr09,BrIaTr10,BrIaTr11,BrIaTr12,BrIaTr12a} (see also \cite{BuBu12} for an application of
low-rank R-K methods to stochastic differential equations).

Furthermore, our aim is not to destroy the overall convergence features of the simplified Newton method \eqref{simnewt}.
Thus, instead of \eqref{simnewt1}, we first recast system \eqref{simnewt} as
\begin{equation}
\label{newt2}
\left( I - h(L\otimes J)\right) \Delta^{(k)} = h \left((A-L)\otimes J\right)\Delta^{(k)} - G(\bfy^{(k)}),
\end{equation}
and then, we retrieve an approximation of the unknown vector $\Delta^{(k)}$  by means of the {\em inner iteration}
\begin{equation}
\label{newt3}
\left( I - h(L\otimes J)\right)\Delta^{(k)}_{\nu+1} = h \left((A-L)\otimes J\right)\Delta^{(k)}_\nu - G(\bfy^{(k)}),
\end{equation}
starting at $\Delta^{(k)}_0=0$. The inner scheme \eqref{newt3} could be iterated to convergence or stopped after a
suitable number, say $r$, of steps. We see that \eqref{simnewt1} corresponds to \eqref{newt3} performed with
 one single inner iteration. Considering that no function evaluations are needed during the implementation
  of \eqref{newt3}, we aim to perform the minimum number $r$ of inner iterations that does not alter the convergence
   rate of the outer iteration \eqref{newt1}.

The convergence properties of the purely linear scheme \eqref{newt3} continue to be described by the
amplification matrix $M(q)$ defined at \eqref{linerr}. In fact, its iteration matrix is
$$h( I - h(L\otimes J))^{-1}((A-L)\otimes J),$$
 which reduces to $M(q)$ for the individual components corresponding to the eigenvalues
 $\lambda$ of $J$. An advantage of the change of variable we propose is that a fast convergence rate is
 guaranteed at the very first steps of the process, and we will show that, in many practical situations,
 choosing $\nu\le s$ produces very good results (see Table \ref{parametri1}).

The paper is organized as follows. The low-rank formulation of Gauss Radau IIA methods is presented
in Section~\ref{two}, while the splitting procedure is defined in Section~\ref{three}.
Its convergence analysis and some comparisons with similar splitting procedures are reported in Section~\ref{four}.
Section~\ref{five} is  devoted to some numerical tests with the fortran 77 code RADAU5 \cite{HaWa91,testset},
modified according to the presented procedure. Finally, a few conclusions are reported in Section~\ref{six},
along with future directions of investigations.

\section{Augmented low-rank implementation of Radau IIA methods}\label{two}

The discrete problem generated by the application of an $s$-stage ($s\ge2$) Radau IIA method to problem
(\ref{ivp}) may be cast in vector form, by using the W-transform \cite{HaWa91}, as:
\begin{equation}\label{radau}
\bfy = \bfe\otimes y_0 + h \PP X_s \PP^{-1}\otimes I\, f(\bfy),
\end{equation}
where $\bfe$, $\bfy$ and $f(\bfy)$ are defined at \eqref{vectors}, while the matrices $\PP$ and $X_{s}$ are defined as
\begin{equation}\label{PX}
\PP = \pmatrix{ccc}
P_0(c_1) & \dots & P_{s-1}(c_1)\\
\vdots   &           &\vdots\\
P_0(c_s) & \dots & P_{s-1}(c_s)\endpmatrix, \quad
X_s = \pmatrix{ccccc}
\frac{1}2 & -\xi_1\\
\xi_1      &0          &\ddots\\
              &\ddots  &\ddots &-\xi_{s-1}\\
              &            &\xi_{s-2} &0 &-\xi_{s-1}\\
              &            &               &\xi_{s-1} &\bb_s\endpmatrix,\end{equation}
with $\{P_j\}$ the shifted and normalized Legendre polynomials on the interval $[0,1]$,
$$\int_0^1P_i(x)P_j(x)\dd x = \dd_{ij}, \qquad i,j\ge0,$$  and
$$\xi_i = \frac{1}{2\sqrt{4i^2-1}}, \quad i=1,\dots,s-1, \qquad \bb_s =  \frac{1}{4s-2}.$$
Clearly, $h$ is the step size and the abscissae $\{c_1,\dots,c_s\}$ are the Gauss-Radau nodes in $[0,1]$. In particular,
$c_s=1$, so that $y_s$ is the approximation to the true solution at the time $t_{1}=t_{0}+h$.

We now derive an augmented low-rank Runge-Kutta method, which is equivalent to (\ref{radau}),
 by following an approach similar to that devised in \cite{BrIaTr09} to introduce {\em Hamiltonian boundary value methods}
 (HBVMs), a class of energy-preserving R-K methods. In more detail, we choose an auxiliary set of distinct abscissae,
\begin{equation}\label{hc}
0<\hc_1<\dots<\hc_s=1,
\end{equation}
and define the following change of variables involving the internal stages $y_{i}$:
\begin{equation}\label{radau1_2}
\hbfy = \hPP \PP^{-1}\otimes I\,\bfy,
\end{equation}
with
$$\hbfy =\pmatrix{c} \hy_1\\ \vdots \\ \hy_s\endpmatrix,  \qquad \hPP = \pmatrix{ccc}
P_0(\hc_1) & \dots & P_{s-1}(\hc_1)\\
\vdots   &           &\vdots\\
P_0(\hc_s) & \dots & P_{s-1}(\hc_s)\endpmatrix.$$
The vectors $\{\hy_i\}$, $i=1,\dots,s$,  called {\em auxiliary stages},\footnote{They are called {\em silent stages}
in the HBVMs terminology, since their presence does not alter the complexity of the resulting nonlinear system.
 Similarly, the abscissae (\ref{hc}) are called {\em silent abscissae}.} are nothing but the values
 at the abscissae (\ref{hc}) of the polynomial interpolating  the internal stages $\{y_i\}$.
 Substituting \eqref{radau1_2} into \eqref{R-K} yields the new nonlinear system in the unknown $\hbfy$
 (notice that $\hPP \PP^{-1}\bfe = \bfe$):
\begin{equation}
\label{radau1_1}
\hG(\hbfy) \equiv \hbfy-\bfe\otimes y_0-h\hPP X_s\PP^{-1} \otimes I\, f\left( \PP\hPP^{-1}\otimes I\,\hbfy\right) = \bf0.
\end{equation}
Of course, after computing $\hbfy$, the solution must be advanced in the standard manner, that is by means of the last component, $y_s$, of the original stage vector $\bfy$. However notice that $\hat c_s= c_s \Rightarrow \hat y_s=y_s$, so that this step of the procedure is costless.

In the next section, we show that the auxiliary abscissae (\ref{hc}) can be chosen so that the solution of the
corresponding simplified Newton iteration (see (\ref{new}) below) is more efficient than solving (\ref{simnewt}).
We end this section by  noticing that system \eqref{radau1_1} is actually identified by a R-K method with rank deficient
coefficient matrix.
\begin{theo} The method (\ref{radau1_2})-(\ref{radau1_1}) can be cast as a Runge-Kutta method with $2s$-stages,
defined by the following Butcher tableau:
$$\begin{array}{c|cc}
\hbfc   & O & \hPP X_s \PP^{-1}\\
\bfc     & O & \PP X_s \PP^{-1}\\
\hline
           & {\bf0}^T &\bfb^T\end{array}$$
where $\bfc$, $\hbfc$ are the vectors with the Radau abscissae and the auxiliary abscissae (\ref{hc}), respectively,
and $\bfb$ contains the weights of the Radau quadrature.
\end{theo}

\section{The splitting procedure}\label{three}
The  simplified Newton iteration (see \eqref{simnewt}) applied to \eqref{radau1_1} reads
\begin{equation}
\label{new}
\begin{array}{l}
\left(I -h\hPP X_s\hPP^{-1}\otimes J\right)\hDelta^{(k)} = -\hG(\hbfy^k), \\
\hbfy^{(k+1)}=\hbfy^{(k)}+\hDelta^{(k)}.
\end{array}
\end{equation}
As we can see, its structure is precisely the same as that we would obtain by applying the simplified Newton
iteration directly to the original system \eqref{radau}, with the only difference that the matrix $\hPP$ in (\ref{new})
should be replaced by $\PP$.

As was emphasized in the introduction, to simplify the structure of systems such as \eqref{new},
van der Houwen and de Swart \cite{HoSw97,HoSw97a} proposed to replace the matrix $(\PP X_s \PP^{-1})$
in (\ref{radau}) with the lower triangular matrix $L$ arising from its  Crout factorization.
The advantage is that, in such a case, to perform the iteration,  one has to factorize $s$
matrices having the same size $m$ as that of the continuous problem with a noticeable saving of work.
They show that on parallel computers this approach gives very interesting speedups over more standard
approaches based on the use of the $LU$ factorization. This is symptomatic of the fact that $LU$ factorizations
generally give a relevant contribution to the overall execution time of a given code.

Similarly, here we want to take advantage from both the Crout factorization of $(\hPP X_s \hPP^{-1})$
appearing in (\ref{new}) and the freedom of choosing the auxiliary abscissae $\{\hc_{i}\}$,  to devise an
iteration scheme that only requires a single $LU$ factorization of a system of dimension $m$ which is,
therefore, suitable for sequential programming. Differently from \cite{HoSw97}, we continue to adopt
the iteration \eqref{new} (outer iteration) and retrieve an approximation to $\hDelta^{(k)}$ via the linear inner iteration
\begin{equation}\label{splitnew}
\left(I-h\hL\otimes J\right) \hDelta^{(k)}_{\nu+1} =
h\left((\hPP X_s \hPP^{-1}-\hL)\otimes J\right)\,\hDelta^{(k)}_{\nu} -\hG(\hbfy^{(k)}),\qquad \nu=0,1,\dots,
\end{equation}
where
\begin{equation}\label{crout}
 \hPP X_s \hPP^{-1} = \hL\hU,
\end{equation}
with $\hL$ lower triangular and $\hU$ upper triangular with unit diagonal entries. Our purpose is to choose the
auxiliary abscissae (\ref{hc}) so that all the diagonal entries of $\hL$ are equal to each other, i.e.,
\begin{equation}\label{allequal}
(\hL)_{jj} = ~^s\sqrt{\det(X_s)}, \qquad j=1,\dots,s.\end{equation}
In so doing, one has to factor only {\em one} $m\times m$ matrix, to carry out the inner iteration (\ref{splitnew}).
 Concerning the diagonal entry in (\ref{allequal}), the following result can be proved by induction.

\begin{theo} Let $X_s$ be defined according to (\ref{PX}) and let
$$\eta = 1+2\lfloor\frac{s}2\rfloor-s \equiv \left\{\begin{array}{cr} 1, &\mbox{\rm if $s$ is even,}\\
0, &\mbox{\rm otherwise,}\end{array}\right.$$ with $\lfloor\cdot\rfloor$ the floor function. Then
\begin{equation}
\label{detXs}
\det(X_s) = \frac{2^{1-s}}{\prod_{i=2-\eta}^{2\lfloor\frac{s}2\rfloor-\eta}(4i^2-1)}.
\end{equation}
Consequently, from (\ref{allequal}) one has:
\begin{equation}\label{diagos}
d_s := (\hL)_{jj}  = \frac{2^{\frac{1}s-1}}{\left(\prod_{i=2-\eta}^{2\lfloor\frac{s}2\rfloor-\eta}(4i^2-1)\right)^{\frac{1}s}},
\qquad j=1,\dots,s.
\end{equation}
\end{theo}
In Table~\ref{cds} we list the auxiliary abscissae $\{\hc_i\}_{i=1,\dots,s}$ and the diagonal entries $d_s$,
given by (\ref{diagos}), for the Radau IIA methods with  $s=2,\dots,5$ stages.
Notice that, having set $\hc_{s}=1$, the free parameters are $s-1$, namely  $\hc_{i}$, $i=1,\dots,s-1$.
We have formally derived the expression of the first $s-1$  diagonal entries of the matrix $\hL$ as a function
of these unknowns, $(\hL)_{jj}\equiv \ell_{j}(\hc_{1},\dots,\hc_{s})$,  and then we have solved the $(s-1)$-dimensional
system $\ell_{j}(\hc_{1},\dots,\hc_{s}) = d_{s}$, $j=1,\dots,s-1$, with the aid of the symbolic computation software Maple.
From \eqref{detXs} it is clear that the last diagonal element of $\hL$ will be automatically equal to $d_{s}$, too.

As was observed in \cite{Co91} in the context of singly implicit R-K methods, the implementation of
a formula such as  \eqref{splitnew} consists of a block-forward substitution which requires the computation
of ~$(T\otimes J)\hDelta^{(k)}_{\nu+1}$,~ with $$T=\hL-d_{s}I$$ (i.e., the strictly lower triangular part of matrix $\hL$),
 at a cost of $O(s^{2}m+m^{2}s)$ operations.  The $O(m^{2}s)$ term, as well as the $m^{2}$
 multiplications for computing ~$(hd_{s})J$~ before the factorization of the matrix ~$I-hd_{s}J$,~ may be
 eliminated by multiplying both sides of \eqref{splitnew} by $$h^{-1}\hL^{-1}\otimes I.$$
 Considering that $$\hL^{-1}=d_{s}^{-1}I-S,$$ with $S$ strictly lower triangular, system \eqref{splitnew} then takes
 the form
\begin{equation}
\label{splitnew1}
\left(\frac{1}{hd_{s}}I-I\otimes J \right) \hDelta^{(k)}_{\nu+1} \\ = \frac{1}{h}(S \otimes I)\hDelta^{(k)}_{\nu+1}+
  \left(C\otimes J\right)\,\hDelta^{(k)}_{\nu} +R^{(k)},
\end{equation}
where $$C=\hL^{-1}(\hPP X_s \hPP^{-1}-\hL)=\hU-I \qquad\mbox{and}\qquad R^{(k)}=
-\frac{1}{h}(\hL^{-1}\otimes I)\hG(\hbfy^{(k)}).$$
Notice that, since $C$ is strictly upper triangular, the multiplication of $J$ by the
 first block-component of $\hDelta^{(k)}_{\nu}$ may be skipped. But we can go another step beyond and
  completely eliminate any $O(m^{2})$ term in the computation of the term $$(C\otimes J)\,\hDelta^{(k)}_{\nu}$$
  at right-hand side of \eqref{splitnew1}. This is true at the very first step, since, by definition, $$\hDelta^{(k)}_{0}=0.$$
 Let us set $$\bfw_{\nu}:=\left(C\otimes J\right)\,\hDelta^{(k)}_{\nu} +R^{(k)},$$ which is part of the right-hand side of \eqref{splitnew1}. Thus $\bfw_0=R^{(k)}$ and the first step of \eqref{splitnew1} is  equivalent to the system
\begin{equation}
\label{firststep}
\left[(hd_{s})^{-1}I-I\otimes J\right] \hDelta^{(k)}_{1}=h^{-1}(S \otimes I)\hDelta^{(k)}_{1}+\bfw_{0}.
\end{equation}
After solving for the unknown $\hDelta^{(k)}_{1}$, we set $\bfv_1$ equal to the right-hand side of \eqref{firststep}, which can be exploited to compute the term
$$(I\otimes J) \hDelta^{(k)}_{1}=(hd_{s})^{-1}\hDelta^{(k)}_{1}-\bfv_{1},$$ at a cost of $O(ms)$ operations. It follows that
\begin{equation}
\label{fuckJ}
(C\otimes J)\,\hDelta^{(k)}_{1} = (C\otimes I)\left[(I \otimes J)\,\hDelta^{(k)}_{1}\right] =
(C \otimes I)\left[(hd_{s})^{-1}\hDelta^{(k)}_{1}-\bfv_{1}\right],
\end{equation}
and thus $\bfw_{1}=\left(C\otimes J\right)\,\hDelta^{(k)}_{1} +R^{(k)}$ may be computed with $O(s^{2}m)$ floating point operations.
This trick may be repeated at the subsequent steps, thus resulting in the following algorithm:
$$
\begin{array}{rcl}
\bfw_0 &:=& R^{(k)}\\
\mbox{do~} \nu&=&0,1,\dots\\
&&\mbox{solve:~} \left[ (hd_s)^{-1}I -I\otimes J\right]\hDelta_{\nu+1}^{(k)} = h^{-1}(S \otimes I)\hDelta^{(k)}_{\nu+1} + \bfw_\nu\\
&&\bfv_{\nu+1} ~:=~ h^{-1}(S \otimes I)\hDelta^{(k)}_{\nu+1} + \bfw_\nu\\
&&\bfw_{\nu+1} ~:=~ (C\otimes I )\left[ (hd_s)^{-1}\hDelta_{\nu+1}^{(k)} -\bfv_{\nu+1}\right] + R^{(k)}\\
\mbox{end}&\mbox{do}&
\end{array}
$$
Notice that $\bfv_{\nu+1}$ is just the right-hand side of the preceding linear system and thus it is freely available as soon as the system has been solved.
\begin{table}[t]
\caption{Auxiliary abscissae for the $s$-stage Radau method, $s=2,\dots,5$,
and the diagonal entry $\ell_{s}$ (see \ref{diagos}) of the corresponding factor $\hL$.}
\label{cds}
\vspace{2mm}
\centerline{\begin{tabular}{|c|c|}
\hline
\multicolumn{2}{|c|}{$s=2$}\\
\hline
$\hc_1$ & $(6-\sqrt{6})/(6+2\sqrt{6})$ \\
$\hc_2$ & 1\\
$d_2$ & 0.40824829046386301636621401245098\\
\hline
\multicolumn{2}{|c|}{$s=3$}\\
\hline
$\hc_1$ & 0.18589230221764097222357873465176\\
$\hc_2$ & 0.50022434784008286059148415923632\\
$\hc_3$ & 1\\
$d_3$ & 0.25543647746451770219954184281099\\
\hline
\multicolumn{2}{|c|}{$s=4$}\\
\hline
$\hc_1$ & 0.12661575733255931078112184952036\\
$\hc_2$ & 0.34154548143311325099490740728171\\
$\hc_3$ & 0.56937072098419698874387077046544\\
$\hc_4$ & 1\\
$d_4$ & 0.18575057999133599176307088298897\\
\hline
\multicolumn{2}{|c|}{$s=5$}\\
\hline
$\hc_1$ & 0.09527975140867214336447374571157\\
$\hc_2$ & 0.28143874673988994521203045137949\\
$\hc_3$ & 0.38152142820340929736570124768463\\
$\hc_4$ & 0.60680555490108389442461323421422\\
$\hc_5$ &1\\
$d_5$ & 0.14591154019899779261811749554182\\
\hline
\end{tabular}}
\end{table}

\section{Convergence analysis and comparisons}\label{four}

In this section we briefly analyze the splitting procedure (\ref{splitnew}). This will be done according to
the linear analysis of convergence in \cite{HoSw97} (see also \cite{BrMa09}). In such a case, problem (\ref{ivp})
becomes the celebrated test equation
\begin{equation}\label{test}
y' = \lambda y, \qquad y(t_0)=y_0.
\end{equation}
By setting, as usual, $q=h\lambda$, one then obtains that the error equation associated with (\ref{splitnew}) is given by
\begin{equation}\label{erreq}
e_{\nu+1} = \hat{M}(q)e_{\nu}, \qquad \hat{M}(q):= q(I-q\hat{L})^{-1}\hat{L}(\hat{U}-I), \qquad \nu=0,1,\dots,
\end{equation}
where we have set $e_{\nu}=\Delta^{(k)}_\nu-\Delta^{(k)}$, that is the error vector at step $\nu$
(we neglect, for sake of simplicity, the index $k$ of the {\em outer} iteration) and $\hat{M}(q)$ is the
 iteration matrix induced by the splitting procedure. This latter will converge if and only if its spectral radius,
$$\rho(q) := \rho(\hat{M}(q)),$$
is less than 1. The {\em region of convergence} of the iteration is then defined as
$$\DD = \left\{ q\in\CC\,:\, \rho(q)<1\right\}.$$
The iteration is said to be {\em $A$-convergent} if $\CC^-\subseteq\DD$. If, in addition, the {\em stiff amplification factor},
$$\rho^\infty := \lim_{q\rightarrow\infty} \rho(q),$$
is null, then the iteration is said to be {\em $L$-convergent}. Clearly, $A$-convergent iterations are
appropriate for $A$-stable methods, and $L$-convergent iterations are appropriate for $L$-stable methods.
In our case, since
\begin{equation}\label{Minfty}
\hat{M}(q) \rightarrow (\hat{U}-I), \qquad q\rightarrow\infty,
\end{equation}
which is a nilpotent matrix of index $s$, the iteration is $L$-convergent if and only if it is $A$-convergent.
Since the iteration is well defined for all $q\in\CC^-$ (due to the fact that the diagonal entry of $\hat{L}$, $d_s$, is positive)
 and $\rho(0)=0,$ $A$-conergence, in turn, is equivalent to require that the {\em maximum amplification factor},
$$\rho^* = \max_{x\in\RR} \rho(ix),$$
is not larger than 1. Another useful parameter is the {\em nonstiff amplification factor},
\begin{equation}\label{M0}
\tilde\rho = \rho(  \hat{L}(\hat{U}-I) ),\end{equation}
that governs the convergence for small values of $q$ since
$$\rho(q) \approx \tilde\rho q, \qquad \mbox{for} \, q\approx 0.$$
Clearly, the smaller $\rho^*$ and $\tilde\rho$, the better the convergence properties of the iteration.
In Table~\ref{parametri} we list the nonstiff amplification factors and the
maximum amplification factors for the following $L$-convergent iterations applied to the $s$-stage Radau IIA methods:
\begin{itemize}
\item[(i)] the iteration obtained by the original triangular splitting in \cite{HoSw97};
\item[(ii)] the iteration obtained by the modified triangular splitting in \cite{AmBr97};
\item[(iii)] the {\em blended} iteration obtained by the {\em blended implementation} of the methods, as defined in \cite{BrMa02};
\item[(iv)] the iteration defined by (\ref{splitnew}).

\end{itemize}
We recall that the scheme (i) (first column) requires $s$ real factorizations per iteration, whereas (ii)--(iv) only need one factorization per iteration.  From the parameters listed in the table, one concludes that the proposed splitting procedure is the most effective among all the considered ones.

It is worth mentioning that the above amplification factors are defined in terms of the eigenvalues of the involved matrices.
Therefore, they are significant if a large number of inner iterations are performed or if the initial guess is accurate enough. In the computational practice, the number of inner iteration is usually small, so that it is also useful to check
the so called {\em averaged amplification factors} over $\nu$ iterations, defined as follows (see (\ref{M0})
and (\ref{Minfty})):
$$\tilde{\rho}_\nu = ~^\nu\sqrt{ \left\| \left[\hat{L}(\hat{U}-I)\right]^\nu\right\|}, \qquad
\rho_\nu^* = \max_{x\in\RR} ~^\nu\sqrt{\|M(ix)^\nu\|}, \qquad \rho_\nu^\infty = ~^\nu\sqrt{ \left\| (\hat{U}-I)^\nu\right\|}.$$
Clearly, $$\rho_\nu^\infty = 0, \qquad \forall \nu\ge s,$$ since matrix $\hat{U}-I$ is nilpotent of index $s$.
Moreover,  $$\tilde\rho_\nu\rightarrow \tilde\rho, \qquad \rho_\nu^*\rightarrow\rho^*,
\qquad \mbox{as}\qquad \nu\rightarrow\infty.$$
For this reason, in Table~\ref{parametri1} we compare the asymptotic parameters
$\tilde\rho$ and $\rho^*$ (columns 2 and 3) with the averaged ones over $s$ iterations
(columns 4 and 5), for $s=2,\dots,5$. As one can see, the iterations are still $L$-convergent
after $s$ iterations (the norm $\|\cdot\|_\infty$ has been considered).
In the last three columns of the table, we list the amplification factors after just 1 inner iteration: in such a case,
the iterations are no more $L$-convergent, though still $A$-convergent, up to $s=4$.

\begin{table}[t]
\begin{center}
\caption{Amplification factors for the triangular splitting in \cite{HoSw97},
the modified triangular splitting in \cite{AmBr97}, the {\em blended} iteration in \cite{BrMa02},
 and the splitting (\ref{splitnew}), for the $s$-stage Radau IIA methods.}
\label{parametri}

\vspace{2mm}
\begin{tabular}{|c|cc|cc|cc|cc|}
\hline
       &\multicolumn{2}{c|}{(i): triangular}
       &\multicolumn{2}{c|}{(ii): triangular}
       &\multicolumn{2}{c|}{(iii): {\em blended}}
       &\multicolumn{2}{c|}{(iv): triangular}\\
       &\multicolumn{2}{c|}{splitting in \cite{HoSw97}}
       &\multicolumn{2}{c|}{splitting in \cite{AmBr97}}
       &\multicolumn{2}{c|}{iteration in \cite{BrMa02}}
       &\multicolumn{2}{c|}{splitting (\ref{splitnew})}\\
\hline
$s$ & $\tilde{\rho}$ & $\rho^*$  & $\tilde{\rho}$ & $\rho^*$ & $\tilde{\rho}$ & $\rho^*$ & $\tilde{\rho}$ & $\rho^*$\\
\hline
2   &  0.1500 & 0.1837 & 0.1498 & 0.1835  &  0.1498 &  0.1835 & 0.1498 &   0.1835\\
3  &  0.1853 & 0.3726 & 0.1375 & 0.3138  &   0.1674 &  0.3398 &  0.1333 &   0.3134  \\
4  &  0.1728 & 0.5064 & 0.1236 & 0.4137 &   0.1535  &  0.4416 &   0.1174 &   0.3826  \\
5  &  0.1496 & 0.6103 & 0.1090 & 0.4949 &   0.1367  &  0.5123 &  0.0787 &   0.3963  \\
\hline
\end{tabular}\\

\bigskip
\caption{Amplification factors, and averaged amplification factors after $s$ inner iterations and 1 inner iteration,
for the triangular splitting  (\ref{splitnew}), for the $s$-stage Radau IIA methods.}
\label{parametri1}

\vspace{2mm}
\begin{tabular}{|c|cc|cc|ccc|}
\hline
$s$ & $\tilde{\rho}$ & $\rho^*$  & $\tilde{\rho}_s$ & $\rho_s^*$ & $\tilde{\rho}_1$ & $\rho_1^*$ & $\rho_1^\infty$\\
\hline
2   & 0.1498 &   0.1835  & 0.1498  &  0.1835 & 0.1498  &  0.2020  &  0.2020\\
3  &  0.1333 &   0.3134  & 0.1407  &  0.3378 & 0.1513  &  0.3984  &  0.3440\\
4  &  0.1174 &   0.3826  & 0.1316  &  0.4363 & 0.2169  &  0.6643  &  0.5172\\
5  &  0.0787 &   0.3963  & 0.1200  &  0.5841 & 0.2959  &  1.1141  &   0.9945\\
\hline
\end{tabular}
\end{center}
\end{table}

\section{Numerical Tests}\label{five}
In this section, we report a few results on three stiff problems taken from the {\em Test Set for IVP Solvers} \cite{testset}:
\begin{itemize}
\item  {\em Elastic Beam problem}, of dimension $m=80$;

\item {\em Emep problem}, of dimension $m=66$;

\item {\em Ring Modulator problem}, of dimension $m=15$.
\end{itemize}
All problems have been solved by using the RADAU5 code \cite{HaWa91,testset}
and a suitable modification of it which implements the splitting procedure
with a fixed number of inner iterations, namely $\nu=1,2,3.$

Clearly, further improvements could be obtained by dynamically varying the number
of inner iterations as well as by implementing a suitable strategy, well tuned for the new iterative procedure,
 to decide whether the evaluation of the Jacobian can be avoided.  In absence of such refinements, in order to verify the effectiveness of the proposed approach, we have forced the evaluation of the Jacobian after every accepted step by setting in input {\tt work(3)=-1D0}. As a consequence,  the factorization of the involved matrices is computed at each integration step.

All the experiments have been done on a PC with an Intel Core2 Quad Q9400 @ 2.66GHz processor under Linux by using the GNU Fortran
compiler {\tt gfortran} with optimization flag {\tt -Ofast}.

The following input tolerances for the relative ($rtol$) and absolute ($atol$) errors and  initial stepsizes ($h_0$)
have been used:
\begin{itemize}
\item Elastic Beam problem: ~$rtol = atol = h_0 = 10^{-4 -i/4}$,~ $i=0,\dots,16$;
\item Emep problem: ~ $rtol = 10^{-4 -i/4}$,~ $i=0,\dots,28,$~$atol=1$ and $h_0 = 10^{-7}$;
\item Ring Modulator problem: ~$rtol = atol = h_0 = 10^{-7-i/4}$,~ $i=0,\dots,20$.
\end{itemize}

\begin{figure}
\centerline{\includegraphics[width=14cm,height=8cm]{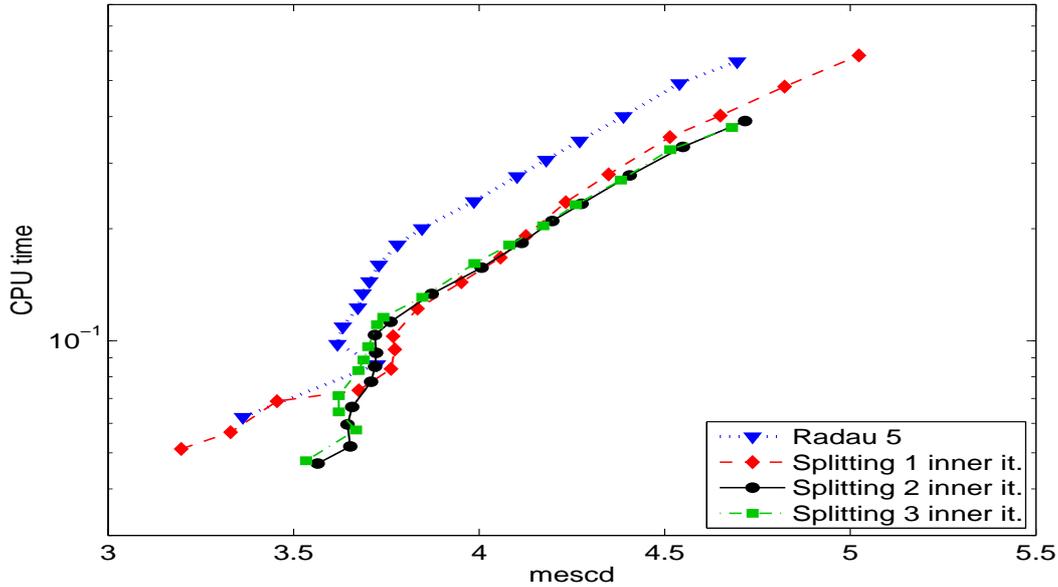}}
\caption{Work precision diagram for the Elastic Beam problem.}
\label{beam}
\end{figure}

\begin{figure}
\centerline{\includegraphics[width=14cm,height=8cm]{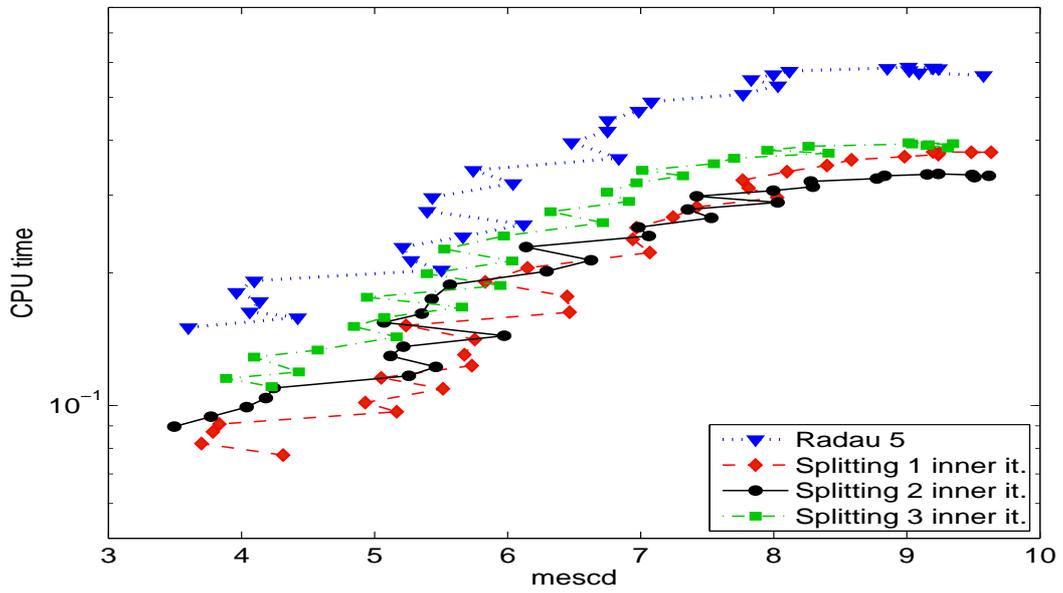}}
\caption{Work precision diagram for the Emep problem.}
\label{emep}
\end{figure}

\begin{figure}
\centerline{\includegraphics[width=14cm,height=8cm]{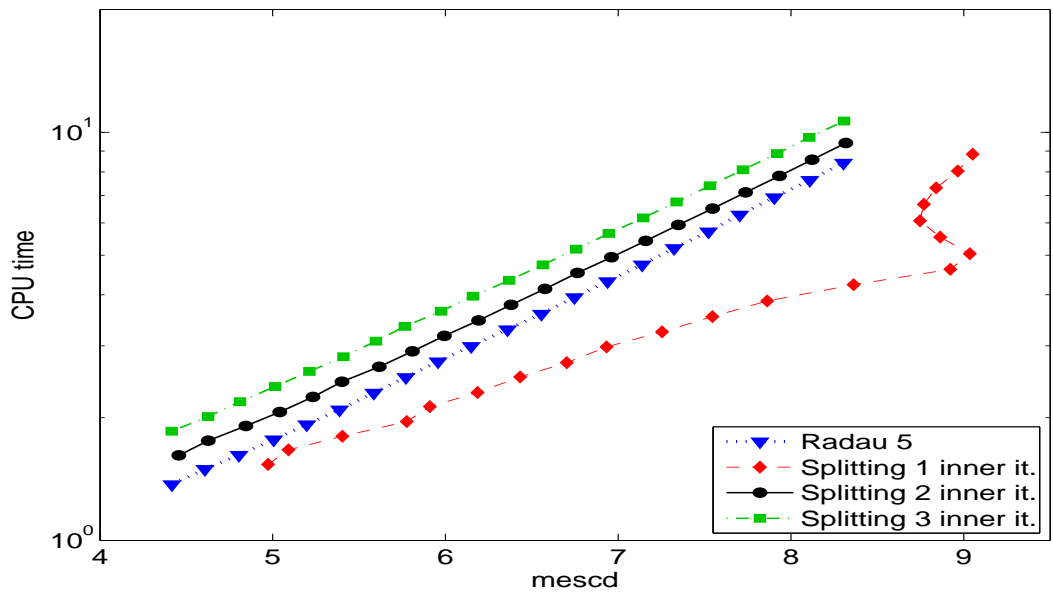}}
\caption{Work precision diagram for the Ring Modulator problem.}
\label{ring}
\end{figure}

Figures~\ref{beam}, \ref{emep}, and \ref{ring} show the obtained results as {\em work-precision diagrams},
where the CPU-time in seconds is plotted versus accuracy, measured as {\em mixed-error significant correct digits (mescd)},
defined as
$$-\log_{10} \max_{i=1,\dots,m} \left\{|e_i|/(1+|y_i|)\right\},$$
$e_i$ being the error in the $i$th entry of the solution at the end of the trajectory, and $y_i$ the corresponding
reference value (which is known, for all problems in the {\em Test Set}).

For the first two problems, the work-precision diagrams suggest that
the splitting version of the RADAU5 code is more efficient than the original one, even starting with 1 inner iteration.
Moreover, in Tables~\ref{radaut}--\ref{radau3} we list a few statistics for the Elastic Beam problem, from which one
deduces that, by using 2--3 inner iterations, the number of steps is approximately the same as the original code: in other words, the convergence rate of the outer iteration is preserved.

For the last problem (Ring Modulator), which has a much smaller size,
the splitting with $2$ and $3$ inner iterations is less efficient than the original RADAU5 code.
Nevertheless, when using a single inner iteration  the algorithm uses a larger number of steps (8-10\% more), as is shown in Tables~\ref{ring0} and \ref{ring1}, resulting into a much more accurate solution. In our understanding, this behaviour may be explained by considering that computing the vector field $f(t,\bfy)$ of this problem is extremely cheap and hence accuracy is more conveniently obtained by acting on the number of function evaluations rather than on the number of inner iterations.

\begin{table}
\begin{center}
\caption{Statistics for the Elastic Beam problem, RADAU5.}
\label{radaut}
\begin{tabular}{|rrrrrrrr|}
\hline
rtol &mescd &steps &accept &feval &jeval &LU &CPU-time\\
\hline
1.00E-04  & 3.36 &55 &49 &380 &49 &55 &6.24E-02\\
1.00E-05  & 3.67 &112 &95 &764 &95 &112 &1.23E-01\\
1.00E-06  & 3.78 &162 &146 &1103 &146 &162 &1.81E-01\\
1.00E-07  & 4.18 &275 &251 &1853 &251 &275 &3.06E-01\\
1.00E-08  & 4.69 &507 &459 &3417 &459 &507 &5.62E-01\\
\hline
\end{tabular}\\~\\
\caption{Statistics for the Elastic Beam problem, RADAU5, split 1.}
\label{radau1}
\begin{tabular}{|rrrrrrrr|}
\hline
rtol &mescd &steps &accept &feval &jeval &LU &CPU-time\\
\hline
1.00E-04  & 3.20 &74 &66 &870 &66 &74 &5.12E-02\\
1.00E-05  & 3.76 &117 &105 &1443 &105 &117 &8.40E-02\\
1.00E-06  &  3.95 &193 &177 &2769 &177 &193 &1.44E-01\\
1.00E-07  & 4.35 &374 &330 &5925 &330 &374 &2.80E-01\\
1.00E-08  & 5.02 &801 &655 &12814 &655 &801 &5.84E-01\\
\hline
\end{tabular}\\~\\
\caption{Statistics for the Elastic Beam problem, RADAU5, split 2.}
\label{radau2}
\begin{tabular}{|rrrrrrrr|}
\hline
rtol &mescd &steps &accept &feval &jeval &LU &CPU-time\\
\hline
1.00E-04  & 3.57 &66 &56 &548 &56 &66 &4.68E-02\\
1.00E-05  & 3.71 &112 &96 &879 &96 &112 &7.76E-02\\
1.00E-06  & 3.76 &152 &144 &1290 &144 & 152 &1.12E-01\\
1.00E-07  & 4.20 &284 &260 &2603 &260 &284 &2.10E-01\\
1.00E-08  & 4.72 &517 &481 &5044 &481 &517 &3.89E-01\\
\hline
\end{tabular}\\~\\
\caption{Statistics for the Elastic Beam problem, RADAU5, split 3.}
\label{radau3}
\begin{tabular}{|rrrrrrrr|}
\hline
rtol &mescd &steps &accept &feval &jeval &LU &CPU-time\\
\hline
1.00E-04  &3.53 &64 &54 &454 &54 & 64 &4.76E-02\\
1.00E-05  &3.67 &115 &96 &810 &96 &115& 8.32E-02\\
1.00E-06  &3.74 &154 &141 &1104 &141 &154 &1.16E-01\\
1.00E-07  &4.17 &273 &249 &1959 &249 &273 &2.04E-01\\
1.00E-08  &4.68 &502 &456 &3654 &456 &502 &3.74E-01\\
\hline
\end{tabular}
\end{center}
\end{table}

\begin{table}
\begin{center}
\caption{Statistics for the Ring Modulator problem, RADAU5.}
\label{ring0}
\begin{tabular}{|rrrrrrrr|}
\hline
rtol &mescd &steps &accept &feval &jeval &LU &CPU-time\\
\hline
1.00E-07  &  4.42 &   98754 &   89346 &   510295 &   89346 &   98754 & 1.37E+00\\
1.00E-08  &  5.20 & 137823 & 128316 &   727506 & 128316 & 137823 & 1.92E+00\\
1.00E-09  &  5.96 & 194463 & 185008 & 1046747 & 185008 & 194463 & 2.74E+00\\
1.00E-10  &  6.75 & 277830 & 268414 & 1525756 & 268414 & 277830 & 3.94E+00\\
1.00E-11  &  7.52 & 399846 & 390508 & 2234881 & 390508 & 399846 & 5.71E+00\\
1.00E-12  &  8.30 & 580535 & 571309 & 3365783 & 571309 & 580535 & 8.42E+00\\
\hline
\end{tabular}\\~\\
\caption{Statistics for the Ring Modulator problem, RADAU5, split 1.}
\label{ring1}
\begin{tabular}{|rrrrrrrr|}
\hline
rtol &mescd &steps &accept &feval &jeval &LU &CPU-time\\
\hline
1.00E-07  &  4.97 & 110376 &   95269 &   958749 &  95269  & 110376 & 1.54E+00\\
1.00E-08  &  5.91 & 152526 & 136231 & 1328822 & 136231 & 152526 & 2.13E+00\\
1.00E-09  &  6.93 & 212686 & 195982 & 1855438 & 195982 & 212686 & 2.98E+00\\
1.00E-10  &  8.36 & 301719 & 283921 & 2635810 & 283921 & 301719 & 4.23E+00\\
1.00E-11  &  8.75 & 432000 & 412643 & 3785978 & 412643 & 432000 & 6.07E+00\\
1.00E-12  &  9.05 & 624708 & 602385 & 5524392 & 602385 & 624708 & 8.84E+00\\
\hline
\end{tabular}
\end{center}
\end{table}

\section{Conclusions}\label{six}
In this paper we have defined a splitting procedure for Radau IIA methods, derived by an {\em augmented low-rank}
formulation of the methods. In such formulation, a set of {\em auxiliary abscissae} are determined such that
 the Crout factorization of a corresponding matrix associated with the method has constant diagonal entries.
 In such a case, the complexity of the iteration is optimal. Moreover, the presented iteration compares favorably
  with all previously defined iterative procedures for the efficient implementation of Radau IIA methods.
  The presented technique can be straightforwardly extended to other classes of implicit Runge-Kutta methods
  (e.g., collocation methods) and this will be the subject of future investigations.

\end{document}